\documentclass[11pt,twoside]{article}
\usepackage{latexsym}
\usepackage{amssymb,amsbsy,amsmath,amsfonts,amssymb,amscd}
\usepackage{epsfig, graphicx,subfigure}
\usepackage[active]{srcltx}
\usepackage{psfrag}
\usepackage{multicol,float}
\usepackage{colortbl}
\newcommand{\commentout}[1]{}

\newcommand{\1}{{\mathchoice {\rm 1\mskip-4mu l} {\rm 1\mskip-4mu l}
{\rm 1\mskip-4.5mu l} {\rm 1\mskip-5mu l}}}

\newcommand {\al} {\alpha}

\newcommand {\lb} {\lambda}
\newcommand {\Chi} {{\bf \raise 2pt \hbox{$\chi$}} }

\newcommand {\U} { {\mathcal U} }
\newcommand {\f}   {\frac}

\newcommand{\dis}{\displaystyle}

\newcommand{\beq}{\begin{equation}}
\newcommand{\beqa} {\begin{array}{rl}}
\newcommand{\eeq}{\end{equation}}
\newcommand{\eeqa}{\end{array}}

\newtheorem{lemma}{Lemma}

\newtheorem{corollary}{Corollary}
\newcommand{\qed}{{ \hfill
                       {\unskip\kern 6pt\penalty 500
                       \raise -2pt\hbox{\vrule\vbox to 6pt{\hrule width 6pt
                       \vfill\hrule}\vrule} \par}   }}
\title{\LARGE The Shape of the Polymerization Rate in the Prion Equation.}

\author{Pierre Gabriel\thanks{Universit\'e Pierre et Marie Curie-Paris 6, UMR 7598 LJLL, BC187, 4, place de Jussieu, F-75252 Paris cedex 5; email: gabriel@ann.jussieu.fr}}

\date{\today}

\begin{document}
\maketitle
\pagestyle{plain}
\pagenumbering{arabic}

\begin{abstract}
We consider a polymerization (fragmentation) model with size-dependent parameters involved in prion proliferation. Using power laws for the different rates of this model, we recover the shape of the polymerization rate using experimental data. The technique used is inspired from \cite{Legname}, where the fragmentation dependence on prion strains is investigated. Our improvement is to use power laws for the rates, whereas \cite{Legname} used a constant polymerization coefficient and linear fragmentation.
\end{abstract}

\


\

\noindent{\bf AMS Class. No.} 35F50, 35Q92, 35R30, 45C05, 45K05, 45M99, 45Q05, 65D10, 92D25

\

\section{Introduction}

Polymerization of prion proteins is a central event in the mechanism of prion diseases. The cells naturally produce the normal form of this protein (PrPc), but there exists a pathogenic form (PrPsc) present in infected cells as polymers. These polymers have the ability to transconform and attach the normal form by a polymerization process not yet very well understood. Mathematical modelling is a promising tool to investigate this mechanism, provided we can find accurate parameters.\\
We consider the following model of prion proliferation (see \cite{CL2,CL1,Greer,LW}):
\begin{equation}
\label{eq:model}
\left\{
\begin{array}{rll}
\dfrac{dV(t)}{dt}&=&\displaystyle \lambda- V(t)\left[\delta+ \int_{0}^{\infty} \tau(x) u(x,t) \; dx\right],
\vspace{.2cm}\\
\dfrac{\partial}{\partial t} u(x,t) &=& \displaystyle - V(t) \frac{\partial}{\partial x} \big(\tau(x) u(x,t)\big) - [\mu_0 + \beta(x)] u(x,t)
+ 2 \int_{x}^{\infty} \beta(y) \kappa (x,y) \, u(y,t) \, dy,
\vspace{.2cm}\\
u(0,t) &=&0,
\vspace{.2cm}\\
u(x,0) &=&u_0(x).
\end{array} \right.
\end{equation}
This equation is the continuous version of the discrete model of \cite{Masel} (see \cite{DGL} for the derivation). System~\eqref{eq:model} is a quite general model able to reproduce biological experimental results (see \cite{Lenuzza}). In this equation, $V(t)$ is the quantity of PrPc which is produced by the cell with a rate $\lb,$ and degraded with the rate $\delta.$ The quantity of aggregates of PrPsc with size $x$ at time $t$ is denoted by $u(x,t).$ Polymers of PrPsc can attach monomers (PrPc) with rate $\tau$ and can split themselves in two smaller polymers with rate $\beta.$ To fit with experimental observations \cite{Silveira}, the polymerization and fragmentation rates should depend on the polymer size $x$ (see \cite{CL2,CL1}). When fragmentation occurs, a polymer of size $y$ produces two pieces of size $x$ and $y-x$ with the probability $\kappa(x,y).$ The natural assumptions (see \cite{MMP1,MMP2,BP,LW}) on the kernel $\kappa$ are then
\beq\label{as:kappasup}Supp\,\kappa\subseteq\{(x,y),\ 0\leq x\leq y\}\eeq
since fragmentation produces smaller aggregates than the initial one,
\beq\label{as:kappa1}\forall y>0,\quad\int\kappa(x,y)\,dx=1\qquad(\kappa(\cdot,y)\ \text{is a probability measure})\eeq
and
\beq\label{as:kappa2}\forall x\leq y,\quad\kappa(x,y)=\kappa(y-x,y)\eeq
because of the indistinguishability of the two pieces of sizes $x$ and $y-x$ after fragmentation. A consequence of the assumptions \eqref{as:kappasup}-\eqref{as:kappa2} is that the mean size of the fragments obtained from a polymer of size $y$ is $y/2:$
\beq\label{eq:kappa3}\forall y>0,\quad\int x\kappa(x,y)\,dx=\f y2.\eeq
In the discrete model of \cite{Masel3,Masel,Masel2}, a critical size $n$ of polymers is considered. Under this size, the polymers are unstable, so the proteins that compose them go back to their normal form PrPc. This effect is taken into account in the continuous model of \cite{Greer,Greer2} with the critical size $x_0.$ As a consequence, there is one more integral term of apparition of monomers in the equation on $V.$ Nevertheless, even if $n>0$ in the discrete model, the value of $x_0$ can be taken equal to $0$ in the continuous one, as justified in \cite{DGL}. Then the fact that small polymers are unstable, and so do not exist, is taken into account by the boundary condition $u(0,t)=0.$ The parameter $\mu_0$ is a death term which will be discussed later.

In the article \cite{Legname}, model \eqref{eq:model} with the specific parameters of \cite{Engler,Greer,Greer2,Masel3,Masel,Masel2} is used. It assumes that the polymerization rate does not depend on polymer size ($\tau(x)\equiv\tau_0$) and the fragmentation is described by $\beta(x)=\beta_0x$ and $\kappa(x,y)=\1_{0\leq x\leq y}\f1y.$ In this case, equation~\eqref{eq:model} can be reduced to a non-linear system of ODEs (see \cite{Engler,Greer,Greer2} or even \cite{Masel3,Masel,Masel2}). This system is used in \cite{Legname} to investigate the parameters involved in the so-called ``prion strain phenomenon''. This phenomenon is the ability of the same encoded protein to encipher a multitude of phenotypic states. Prion strains are characterized by their incubation period (time elapsed between experimental inoculation of the PrPsc infectious agent and clinical onset of the disease) and their lesion profiles (patterns produced by the death of neurones) in the brain \cite{Fraser}. A reason for the strain diversity can be that there exist various conformation states of PrPsc, namely various misfoldings of the prion protein, which lead to various biological and biochemical properties \cite{Morales,Rezaei}. These different conformations could be characterizable by their different stability against denaturation. Indeed, recent studies on prion disease have confirmed that the incubation time is related not only to the quantity of pathogenic protein inoculated ($\int u_0(x)\,dx$) and the level of expression of PrPc ($\lambda$), but also to the resistance of prion strains against denaturation \cite{Legname3} in terms of the concentration of guanide hydrochloride (Gdn-HCl) required to denature $50\%$ of the disease-causing protein. Other studies have highlighted a strong relationship between the stability of the prion protein against denaturation and neuropathological lesion profiles \cite{Legname2,Sigurdson}. A natural question is then to know how the conformation states of the PrPsc protein affect the parameters of model~\eqref{eq:model} and more precisely which of these parameters is the most strain dependent. One of the conclusions of \cite{Legname} is that the key parameter describing strain-dependent replication kinetics is the fibril breakage rate $\beta_0.$ Using the model, some coefficients already defined in \cite{Masel} are calculated from experimental data, and the dependence between them indicates that the fragmentation rate $\beta_0$ is the most adapted parameter to be strain dependent. But, in this framework, the fragmentation dependence is not sufficient to fit the data precisely, and the polymerization rate $\tau_0$ is also assumed to change with the strain. In particular, to keep a one-parameter family of models, the polymerization rate is assumed to depend on the fragmentation one through a power relation:
\beq\label{as:betatau}\exists\phi\ s.t.\ \tau_0={\beta_0}^\phi.\eeq
However, polymerization and fragmentation are two different biochemical processes which have no reason to be linked together.

In this paper we prove that we can get rid of such a relation by considering a larger class of parameters for model \eqref{eq:model}. We take advantage of a general framework developed in \cite{CL2,CL1,LW}, where coefficients can be more general than in \cite{Engler,Greer,Greer2,Masel3,Masel,Masel2,Legname}, and we consider power laws as in \cite{EscoMischler4,EscoMischler3,LW,M1}. Thanks to this extension, we find a shape for the polymerization rate allowing us to fit experimental data with only the fragmentation rate dependent on strain. The experimental data we use are the stability against denaturation $G$ and the growth rate of the quantity of PrPsc ($r$), which is linked to the incubation time. These data are strain dependent, and we can find in \cite{Legname} estimated empirical values for different strains. In Section~\ref{se:asandme}, we present the class of coefficients we consider, and give the necessary assumptions for the mathematical investigations. Most of these assumptions are already made in \cite{Masel, Legname} and, thanks to them, we are able to link $r$ and $G$ to the mathematical model. In Section~\ref{se:results} we give the dependence of these parameters on the polymerization and fragmentation rates. Then we give a relation between them using only the fragmentation dependency on strain and are able to fit data with an accurate shape of the polymerization rate.

\

\section{Assumptions and Method}\label{se:asandme}

To avoid the use of a power relation like \eqref{as:betatau} between the polymerization and the fragmentation, we allow the parameters $\tau$ and $\beta$ to depend on the size $x$ with power laws :
\beq\label{as:homogeneous}\exists\gamma,\nu\ \text{and}\ \beta_0,\tau_0\ \text{s.t.},\quad\forall x>0,\quad\beta(x)=\beta_0x^\gamma\ \text{and}\ \tau(x)=\tau_0x^\nu.\eeq
Concerning the fragmentation kernel $\kappa,$ we assume that there exists a measure $\kappa_0$ on $[0,1]$ such that
\beq\label{as:kappa0}\forall x,y>0,\quad\kappa(x,y)=\f1y\kappa_0\Bigl(\f xy\Bigr).\eeq
Then the natural assumptions on $\kappa$ become
\beq\label{as:kappa02}\int_0^1\kappa_0(z)\,dz=1\quad\text{and}\quad\kappa_0(z)=\kappa_0(1-z).\eeq
We assume that the exponents $\nu,\ \gamma,$ the kernel $\kappa_0$ and the death rate $\mu_0$ are strain independent unlike the ``intensities'' $\beta_0$ and $\tau_0.$ First, we look which of $\beta_0$ and $\tau_0$ is the most characteristic of the strain. Then we find a value of $\nu$ allowing to fit experimental data as well as \cite{Legname} but without assumption \eqref{as:betatau}. We will see in Section~\ref{se:results} that the method naturally gives the parameter $\nu,$ because it turns out that $\gamma$ and $\kappa_0$ do not appear in the final argument.

\

The method proposed in \cite{Masel, Legname} is to investigate the dynamics of the system \eqref{eq:model} linearized at the disease-free equilibrium. Then we can use the eigenelements as suggested in \cite{CL2,CL1}. To make this approximation valid, we make some classical assumptions. First we assume that, in the biological experiments (namely before animal death), the amount of PrPc in cells remains close to the sane value $\bar V=\f\lambda\delta$ (obtained by setting $u\equiv0$ in the first equation of \eqref{eq:model}). Quantity $V$ being assumed to be constant, we consider the eigenelements $r$ and $\U$ of the equation for polymers associated with $\bar V:$
\beq\label{eq:eigen}\left\{\begin{array}{l}
\dis r\U(x)=-\bar V\frac{\partial}{\partial x}\big(\tau(x)\U(x)\big)-[\mu_0+\beta(x)]\U(x)+2\int_{x}^{\infty}\beta(y)\kappa(x,y)\,\U(y)\,dy,
\vspace{.2cm}\\
\dis\U(0)=0,\quad\int\U(x)\,dx=1.
\end{array}\right.\eeq
The principal eigenvalue $r$ represents the growth rate of the solutions, and the eigenvector $\U$ is the corresponding size distribution of polymers. To consider that the size distribution is close to this profile along the experiments, we assume that the initial value $u_0(x)$ is already close to $\U.$ It is a reasonable assumption if contaminations are made from other infected individuals.
Then the solution of \eqref{eq:model} is reduced to
\beq\label{eq:sol} V(t)\simeq\bar V\quad\text{and}\quad u(x,t)\simeq\Bigl[\int u_0(y)dy\Bigr]\,\U(x)e^{r t}.\eeq

To solve our problem, we use the idea that the eigenelements $r$ and $\U$ change with the prion strains (see \cite{CL1, Lenuzza} for instance). Because present experiments do not allow one to measure the size distribution $\U,$ we link it to the stability against denaturation $G$ by the relation
\beq\label{eq:affdep}\exists a>0,\ b\quad\text{s.t.}\quad G=a\bar x+b\eeq
where $\bar x=\int x\,\U(x)\,dx$ is the mean size of the polymers since $\int \U(x)\,dx=1$. Such an affine dependence, suggested in \cite{Legname}, is based on the idea that the longer $\bar x$ is, the larger $G$ will be. In \cite{Legname2,Legname} we can find estimated empirical values of $G$ (as a concentration, in $mol,$ of Gnd-HCL required to denature $50\%$ of PrPsc) and of the growth rate $r$ (computed from measures of incubation times thanks to the assumption of an exponential growth homogeneous in the brain). These values are reported in Table~\ref{tab:data}.

\

\begin{table}[h]
\begin{center}
\begin{tabular}{c||cccccccc}
\hline
Prion strain & 139A & ME7 & BSE & Sc237 & RML & vCJD & Chandler Scrapie & 301 V\\
\hline
r($day^{-1}$) & 0.05 & 0.024 & 0.015 & 0.11 & 0.18 & 0.07 & 0.17 & 0.07\\
\hline
G($M$) & 2 & 2.9 & 2.8 & 1.6 & 1.7 & 1.85 & 2.2 & 2.2\\
\hline
\end{tabular}
\caption{\small\label{tab:data}Estimated empirical values of the growth rate and the stability against denaturation for different prion strains (data from \cite{Legname3,Legname}).}
\end{center}
\end{table}

In the next section, we find a relation between $r$ and $\U$ by investigating their dependence on parameters. Then we show that this relation allows us to fit the data of Table~\ref{tab:data} with an accurate choice of the polymerization exponent $\nu.$

\

\section{Results}\label{se:results}

Proving the existence of solutions to the eigenvalue problem~\eqref{eq:eigen} is not easy and the theory was first developed in \cite{PR, M1}. In \cite{DG}, we could treat singular polymerization and fragmentation rates as investigated here. It is proved in \cite{DG, M1} that solutions exist under the condition
\beq\label{as:gammanu}\gamma+1-\nu>0.\eeq
Moreover, following the idea of \cite{EscoMischler3}, we are able to give the dependence of $r$ and $\U$ on the parameters $V,\ \mu_0,\ \beta_0$ and $\tau_0$ for given $\gamma,\ \nu$ and $\kappa_0$ (unknown but strain independent). Define $r_1$ and $\U_1$ by the eigenequation given for $\beta_0=\bar V\tau_0=1$ and $\mu_0=0$ by
\beq\label{eq:eigeneq1}r_1\U_1(x)=- \frac{\partial}{\partial x} \big(x^\nu \U_1(x)\big) - x^\gamma \U_1(x)
+ 2 \int_{x}^{\infty} y^\gamma \kappa_0\Bigl(\f xy\Bigr) \U_1(y) \, \f{dy}{y}.\eeq
Then $r_1$ and $\U_1$ depend only on $\gamma,\ \nu$ and $\kappa_0,$ and we have the following lemma, which we use later on.

\begin{lemma}\label{lm:autosim}
The eigenelements $r$ and $\U$ solution of \eqref{eq:eigen} satisfy
\beq\label{eq:eigenel2}r=\bigl(V\tau_0\bigr)^{\f{\gamma}{1+\gamma-\nu}}{\beta_0}^{\f{1-\nu}{1+\gamma-\nu}}r_1-\mu_0\quad\text{and}\quad \U(x)=\Bigl(\f{\beta_0}{V\tau_0}\Bigr)^{\f{1}{1+\gamma-\nu}}\U_1\Bigl(\bigl(\f{\beta_0}{V\tau_0}\bigr)^{\f{1}{1+\gamma-\nu}}x\Bigr).\eeq
\end{lemma}

\begin{proof}The proof is based on a homogeneity argument. First we write \eqref{eq:eigen} as
\beq\label{eq:eigeneq2}\frac{r+\mu_0}{V\tau_0}\U(x)=- \frac{\partial}{\partial x} \big(x^\nu \U(x)\big) - \f{\beta_0}{V\tau_0} x^\gamma \U(x)
+ 2 \f{\beta_0}{V\tau_0} \int_{x}^{\infty} y^\gamma \kappa_0\Bigl(\f xy\Bigr) \U(y) \, \f{dy}{y}\eeq
and we set $\al=\f{\beta_0}{V\tau_0}$ and $R=\frac{r+\mu_0}{V\tau_0}.$ The rescaled function $v$ defined by $\U(x)=\al^kv(\al^kx)$ satisfies the equation
\beq\label{eq:eigeneq3}\al^kRv(x)=- \al^{k(2-\nu)}\frac{\partial}{\partial x} \big(x^\nu v(x)\big) - \al^{k+1-k\gamma} x^\gamma v(x)
+ 2 \al^{k+1-k\gamma} \int_{x}^{\infty} y^\gamma \kappa_0\Bigl(\f xy\Bigr) v(y) \, \f{dy}{y}.\eeq
We choose $k=\f{1}{1+\gamma-\nu}$ (which is possible since \eqref{as:gammanu} is satisfied) and obtain
\beq\label{eq:eigeneq4}\al^{\f{\nu-1}{1+\gamma-\nu}}Rv(x)=- \frac{\partial}{\partial x} \big(x^\nu v(x)\big) - x^\gamma v(x)
+ 2 \int_{x}^{\infty} y^\gamma \kappa_0\Bigl(\f xy\Bigr) v(y) \, \f{dy}{y}.\eeq
Comparing \eqref{eq:eigeneq1} and \eqref{eq:eigeneq4}, the uniqueness of eigenelements ensures that
\beq\label{eq:eigenel1}R=\al^{\f{1-\nu}{1+\gamma-\nu}}r_1\quad\text{and}\quad v(x)=\U_1(x)\eeq
because $\int v(x)dx=\int\U(x)dx,$ and the lemma is proved.
\qed
\end{proof}

\

We are now able to give the dependence of $r$ and $\bar x$ on the parameters $\beta_0,\ \tau_0$ and $\mu_0.$ The mean length of the aggregates at time $t$ is defined as the quotient $\dis\frac{\int xu(x,t)\,dx}{\int u(x,t)\,dx}.$ So, using the approximation \eqref{eq:sol}, we obtain that this quantity is constant in time and equal to
\beq\label{def:bars}\bar x=\int x\,\U(x)\,dx.\eeq
Then, thanks to Lemma~\ref{lm:autosim}, we conclude the following corollary.
\begin{corollary}
Under assumptions \eqref{as:homogeneous}-\eqref{as:kappa02} and \eqref{as:gammanu}, there are two constants $r_1$ and $\bar x_1$ (depending only on $\gamma,\ \nu$ and $\kappa_0$) such that
\beq\label{eq:dependency}r=\bigl(V\tau_0\bigr)^{\f{\gamma}{1+\gamma-\nu}}{\beta_0}^{\f{1-\nu}{1+\gamma-\nu}}r_1-\mu_0\qquad\text{and}\qquad\bar x=\Bigl(\f{V\tau_0}{\beta_0}\Bigr)^{\f{1}{1+\gamma-\nu}}\bar x_1.\eeq
\end{corollary}

\begin{proof}
Using Lemma~\ref{lm:autosim}, we immediately have $\dis\bar x=\Bigl(\f{V\tau_0}{\beta_0}\Bigr)^{\f{1}{1+\gamma-\nu}}\bar{x_1}$ with $\dis\bar{x_1}=\int x\,\U_1(x)\,dx.$
\qed
\end{proof}

\

\begin{figure}[H]

\includegraphics[width=\linewidth]{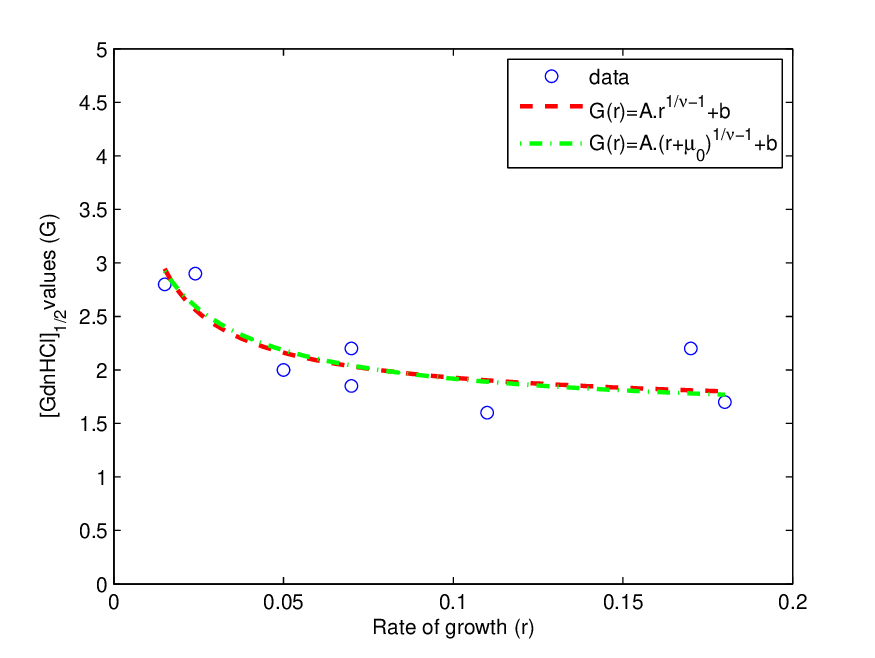}

\caption{\small Experimental values of $G$ and $r$ for different strains are plotted (circles). Then these data are fitted using the reduced model~\eqref{eq:funcdep} with $\mu_0=0$ and $\mu_0$ free.}
\label{fig:fit}

\end{figure}

We have obtained the dependence of $r$ and $\bar x$ on the $a\ priori$ strain-dependent parameters $\beta_0$ and $\tau_0.$ Using assumption~\eqref{eq:affdep}, we deduce the dependence of the stability against denaturation $G$
\beq\label{eq:beta0dep}G=a\Bigl(\f{V\tau_0}{\beta_0}\Bigr)^{\f{1}{1+\gamma-\nu}}\bar{x_1}+b.\eeq
In \cite{Legname}, the exponent $\gamma$ of the fragmentation rate is equal to $1.$ Here we do not impose any value for $\gamma$ but, nevertheless, we consider that it is positive. Indeed the larger the polymers are, the more easily they break themselves. Consequently, a strain-dependency of $\tau_0$ cannot lead to a decreasing correspondence between $r$ and $G,$ as observed in Figure~\ref{fig:fit}, since we have assumed that $\gamma+1-\nu$ and $a$ are positive. Thus it indicates that the key strain-dependent parameter should be $\beta_0.$ With this assumption, we obtain an implicit expression of $G$ as a function of $r:$
\beq\label{eq:funcdep}G=A(r+\mu_0)^\f1{\nu-1}+b\eeq
where $A$ and $b$ are constants (namely strain independent). This model is used to fit the data.

\begin{multicols}{2}

\

\vspace{.9cm}

\begin{table}[H]

\begin{center}
\begin{tabular}{c|cc}
\hline
 & $Ar^\f1{\nu-1}\!+\!b$ & $A(r\!+\!\mu_0)^\f1{\nu-1}\!+\!b$ \\
\hline
$\nu$ & {\bf -0.482} & {\bf 0.316} \\
\hline
$\mu_0$ & 0 & 0.023 \\
\hline
$A$ & 0.083 & 0.01 \\
\hline
$b$ & 1.54 & 1.69 \\
\hline
$R^2$ & 0.7 & 0.72 \\
\hline
\end{tabular}
\vspace{1.56cm}
\caption{\small\label{tab:fit}Best-fitting parameters obtained with the two different models. The correlation coefficient $R^2$ is also reported in the last row.}
\end{center}

\end{table}

\end{multicols}

First we assume as in \cite{Legname} that the death rate $\mu_0$ is equal to zero. Then we obtain a model already considered by \cite{Legname} (Table~3, page~5). In our framework, we can also let $\mu_0$ be free and fit the data with the model~\eqref{eq:funcdep}. As we can see in Table~\ref{tab:fit}, the accuracy of the fitting is a little bit better with this generalization than in the case $\mu_0=0.$ But the most significant difference between the two fittings is the shape of $\tau :$ it is a decreasing function of $x$ in the first case $(\nu<0)$ and it is an increasing function in the second case $(\nu>0).$

\

\section{Conclusion and Perspectives}

Thanks to a larger class of parameters than \cite{Legname}, we are able to fit experimental data by considering only the fragmentation rate to be dependent on strains. We even obtain a more general model and so a more accurate fitting. But, because the parameter profiles are not $a\ priori$ known, we cannot find $\bar x$ from experimental data and have to postulate a particular dependence \eqref{eq:affdep}, whereas this was a consequence of the model in \cite{Legname}. To validate it, new data are necessary, for instance values of $\bar x$ obtained by an experimental size distribution.

A consequence of the method is the prescription of the shape of the polymerization rate from experimental data, and also the death rate $\mu_0.$ Nevertheless, we notice a significant difference between the two recovered exponents $\nu$ obtained by assuming that $\mu_0=0$ or not. In particular, one is negative and the other one positive (see Table~\ref{tab:fit}). The idea of comparing parameters from different strains has been deeply exploited here. But the fact that there are only a few different strains in comparison to the number of fitting parameters is a weakness for precise fitting and precise determination of $\tau.$ Some parameters such as $\mu_0$ or $b$ should be previously determined by another method and from other data. Furthermore, the parameters of model \eqref{eq:model} could have a size dependence more general than that given by power laws. A complete inverse problem is to recover all the size-dependent parameters from experimental data without the restriction to a class of shape. This problem is very complicated, even though first elementary models are treated in \cite{DPZ,PZ,Perasso}. Moreover, it requires more information such as the size distribution $\U$ and not only the mean size $\bar x.$

\vspace{1cm}

{\bf Aknowledgment}

This work has been partly supported by ANR grant project TOPPAZ.

\bibliographystyle{abbrv}
\bibliography{Prion}

\end{document}